\documentclass[preprint,12pt]{elsarticle}
\usepackage{latexsym}
\usepackage[backref=page]{hyperref}
\usepackage{amsmath,amsfonts,amssymb,amscd}
\numberwithin {equation}{section} %
\newtheorem{thm}{Theorem}[section]
\newtheorem{corol}[thm]{Corollary} %
\newtheorem{prop}[thm]{Proposition} %
\newdefinition{rmk}[thm]{Remark} %
\newdefinition{rmks}[thm]{Remarks} %
\newdefinition{example}[thm]{Example}
\newproof{pf}{Proof}
\newcommand{\clos}{\operatorname{cl}}
\newcommand{\ds}{\displaystyle}
\begin{document}

\begin{frontmatter}

\title{Hahn-Banach type extension results for linear operators on asymmetric normed spaces}
\author{S. Cobza\c s}
\address{Babe\c{s}-Bolyai University, Faculty of Mathematics and Computer Science,Department of Mathematics,
Cluj-Napoca, Romania,  E-mail: scobzas@math.ubbcluj.ro}

\begin{abstract}
  We present some results related to Hahn-Banach extension theorem for linear operators  on asymmetric normed spaces. L. Nachbin, Trans. Amer. Math. Soc.  68 (1950), proved
  that a Banach space has the extension  property for linear operators (a property also called injectivity)  if and only if it has the Binary Intersection Property (BIP), meaning that every family of mutually intersecting closed balls has nonempty intersection. Its  analog   for quasi-metric spaces, called mixed BIP,  was considered  by Kemajou et al. Topology Appl. 159   (2012).  The equivalence of mixed BIP to the injectivity of an asymmetric normed space was proved by Conradie et al., Topology Appl. 231 (2017),  derived from some properties of the injective hull of a quasi-metric space.  The aim of the present paper is to give a direct proof of this result   by adapting Nachbin's ideas to the asymmetric case.
\end{abstract}

\end{frontmatter}

\section{Introduction}

A normed space $(Y,\|\cdot\|)$ is said to have the \emph{Hahn-Banach extension property for linear operators} (HB) if for every normed space $(X,\|\cdot\|)$, every closed subspace $Z$ of $X$, every continuous linear operators $T:Z\to Y$ admits a norm-preserving extension to $X$, that is there exists a continuous linear operator $S:X\to Y$ asymmetric normed space
\begin{equation}\label{def.HB}
S|_Z=T\quad\mbox{ and }\quad \|S\|=\|T\|\,.
\end{equation}

Another term used to design this property is that $Y$ is an \emph{injective} normed space.

By the classical Hahn-Banach extension theorem for linear functionals $\mathbb{R}$ and $\mathbb{C}$ are injective.
\begin{rmk}
  Let $X$ be a normed space, $Z$ a subspace of $X$ and $T:Z\to Y$ a continuous linear operator.
 Since the operator $T$ is uniformly continuous, then it has a linear  extension to $\overline Z$ with the same norm, so  that, if the space $Y$ is injective the norm-preserving extension is possible for an arbitrary subspace $Z$ of $X$, not necessarily closed.
\end{rmk}

The normed space $Y$ is said to have the \emph{binary intersection property} (BIP) if  every family $\{B_i\}_{i\in I}$ of closed balls in $Y$  such that $B_i\cap B_j\ne \emptyset,$ for all $i,j\in I,$ has nonempty intersection. It is obvious that the binary intersection property can be defined for metric spaces.
\begin{rmk}\label{re.BIP-compl}
  A metric space having BIP is complete.
\end{rmk}

This follows from the fact that a metric space is complete if and only if  every descending sequence of closed balls with diameters tending to zero has nonempty intersection.

Nachbin \cite{nachbin50} proved the following remarkable result.
\begin{thm}[Nachbin]\label{t1.Nachbin} A Banach space  is injective if and only if it has the BIP
  \end{thm}
It  was  also proved that a Banach space $Y$ is injective if and only if it is isometrically and linearly isomorphic to a space $C(T)$, for an extremally disconnected  compact
Hausdorff space $T$. A topological space $T$ is called extremally  disconnected if the closure of every open subset of $T$ is open. In the real case this was proved by Nachbin \cite{nachbin50} and Goodner \cite{goodner50} under the hypothesis that the unit ball of $Y$ has an extremal point and by Kelley \cite{kelley52} in general. The proof in the complex case was given by Hasumi \cite{hasumi58}. A proof of this result can be found in \cite[Section 4.3]{AK}.

The  study of injectivity in the asymmetric case (meaning quasi-metric spaces and asymmetric normed spaces) was initiated by Salbany \cite{salbany89} and developed in the works of K\"unzi, Otafudu, Agyngi et al, see, for instance, \cite{agyingi17}, \cite{agyingi24}, \cite{otaf-ku-con17}, \cite{otaf-ku-kem12}, and the references quoted therein.

The analog for asymmetric normed spaces of Nachbin's theorem was proved in \cite{otaf-ku-con17}, based of the notion of mixed binary intersection property. The proof was derived from some properties of injective hulls of quasi-metric spaces, developed in \cite{otaf-ku-kem12}. The aim of the present paper is to give a direct proof of this result by adapting the ideas of Nachbin to the asymmetric case.
\section{Quasi-metric spaces and asymmetric normed spaces}
Here we present some notions and results needed for the proof of the main result  (Theorem \ref{t1.asym-Nachbin}).

\subsection{Quasi-metric  spaces}\label{Ss.qpm}

    A {\it quasi-metric} on an arbitrary set $X$ is a mapping $d: X\times X\to
[0,\infty)$ satisfying the following conditions:
\begin{align*}%
\mbox{(QM1)}&\qquad d(x,y)\geq 0, \quad and  \quad d(x,x)=0;\\
\mbox{(QM2)}&\qquad d(x,z)\leq d(x,y)+d(y,z), \\ %
\mbox{(QM3)}&\qquad d(x,y)=d(y,x)=0\Rightarrow x=y,
\end{align*}%
for all $x,y,z\in X.$
The pair $(X,d)$ is called a {\it quasi-metric space}.

The conjugate of the quasi-metric
$d$ is the quasi-metric $\bar d(x,y)=d(y,x),\, x,y\in X.$ The mapping $
d^s(x,y)=\max\{d(x,y),\bar d(x,y)\},\,$ $ x,y\in X,$ is a metric on $X$.

If $(X,d)$ is a quasi-metric space, then for $x\in X$ and $r>0$ we define the balls in $X$ by the formulae %
\begin{align*}%
B_d(x,r)=&\{y\in X : d(x,y)<r\} \; \mbox{-\; the open ball, and }\\ %
B_d[x,r]=&\{y\in X : d(x,y)\leq r\} \; \mbox{-\; the closed ball. } %
\end{align*} %
\begin{rmk}
  The closed ball can  also be defined  for $r=0$, in which case $x\in B_d[x,0],$ but it can contain other points too. If the topology $\tau_d$ is $T_1$  (equivalently, $d(x,y)=0\;\Rightarrow\; y=x$), then $B_d[x,r]=\{x\}.$
\end{rmk}

\textbf{Topological properties}

 The topology $\tau_d$ (or $\tau(d)$) of a quasi-metric space $(X,d)$ can be defined starting from the family
$\mathcal{V}_d(x)$ of neighborhoods  of an arbitrary  point $x\in X$:%
\begin{equation*}
\begin{aligned}
V\in \mathcal{V}_d(x)\;&\iff \; \exists r>0\;\mbox{such that}\; B_d(x,r)\subseteq V\\
                             &\iff \; \exists r'>0\;\mbox{such that}\; B_d[x,r']\subseteq V. %
\end{aligned} %
\end{equation*}

The convergence of a sequence $(x_n)$ to $x$ with respect to $\tau_d,$ called $d$-convergence and
denoted by
$x_n\xrightarrow{d}x,$ can be characterized in the following way %
\begin{equation}\label{char-rho-conv1} %
         x_n\xrightarrow{d}x\;\iff\; d(x,x_n)\to 0. %
\end{equation} %

Also
\begin{equation}\label{char-rho-conv2} %
         x_n\xrightarrow{\bar d}x\;\iff\;\bar d(x,x_n)\to 0\; \iff\; d(x_n,x)\to 0\,.
\end{equation}
and
\begin{equation}\label{char-rho-conv3} \begin{aligned}
 x_n\xrightarrow{d^s}x&\iff\; d^s(x,x_n)\to 0\\&\iff  d(x,x_n)\to 0 \;\mbox{ and }\; d(x_n,x)\to 0\\ &\iff x_n\xrightarrow{d}x\;\mbox{ and }\;  x_n\xrightarrow{\bar d}x\; ,
\end{aligned}\end{equation}

For obvious reasons, some authors call the convergence  \eqref{char-rho-conv1} backward convergence and the convergence \eqref{char-rho-conv2} forward convergence.

The following   topological properties are true for  quasi-metric spaces.
    \begin{prop}[see \cite{Cobzas}]\label{p.top-qsm1}
   If $(X,d)$ is a quasi-metric space, then the following hold.
   \begin{enumerate}
   \item The ball $B_d(x,r)$ is $\tau_d$-open and   the ball $B_d[x,r]$ is
       $\tau_{\bar{d}}$-closed. The ball    $B_d[x,r]$ need not be $\tau_d$-closed.
     \item
   The topology $\tau_d$ is $T_0$ and it is  $T_1$ if and only if
   $d(x,y)>0$ for all  $x\neq y$  in $X$.
      \item   For every fixed $x\in X,$ the mapping $d(x,\cdot):X\to (\mathbb{R},|\cdot|)$ is
   $\tau_d$-usc and $\tau_{\bar d}$-lsc. \\
   For every fixed $y\in X,$ the mapping $d(\cdot,y):X\to (\mathbb{R},|\cdot|)$ is $\tau_d$-lsc and
   $\tau_{\bar d}$-usc.
\end{enumerate}%
     \end{prop} %

As a space equipped with two topologies, $\tau_d$ and   $\tau_{\bar d}\,$, a quasi-metric space can be viewed as a bitopological space in the sense of Kelly \cite{kelly63}.
The following remarks show that imposing too many conditions on a quasi-metric space it becomes metrizable.

\begin{rmk}[\cite{kelly63}] Let $(X,d)$ be  a quasi-metric space. Then
\begin{enumerate}
  \item[\rm (a)]  if  the mapping $d(x,\cdot):X\to (\mathbb{R},|\cdot|)$ is $\tau_d$-continuous for   every $x\in X,$ then the topology $\tau_d$ is regular;
   \item[\rm (b)] if $\tau_d\subseteq \tau_{\bar{d}}$, then the topology $\tau_{\bar{d}}$ is metrizable;
 \item[\rm (c)] if $d(x,\cdot):X\to (\mathbb{R},|\cdot|)$ is $\tau_{\bar{d}}$-continuous for every $x\in X,$ then the topology $\tau_{\bar{d}}$ is metrizable.
\end{enumerate}\end{rmk}

\subsection{Asymmetric normed spaces}

An \emph{asymmetric norm}  on a real vector space $X$ is a function $p:X\to\mathbb{R}_+$ satisfying the conditions:
\begin{itemize}
\item[(AN1)] $p(x)=p(-x)=0\iff x=0$;
\item[(AN2)] $p(\alpha x)=\alpha p(x)$ - positive homogeneity;
\item[(AN3)] $p(x+y)\le p(x)+p(y)$ - the triangle inequality,\end{itemize}
for all $x,y\in X$ and $\alpha\ge 0.$

The pair $(X,p)$ is called an \emph{asymmetric normed space}. The mapping  $\bar p(x)=p(-x),$\, $x\in X,$ is also an asymmetric norm on $X$, called the conjugate of $p$, and $p^s(x)=\max\{p(x),p(-x)\},$\, $x\in X,$ is a norm on $X$. If $(X,p)$ is an asymmetric normed space, then  the mapping $d_p:X\times X\to\mathbb{R}_+$ given by
$d_p(x,y)=p(y-x),\, x,y\in X,$ is a quasi-metric on $X$ and all  topological notions concerning $p$ are considered as those  corresponding to this associated quasi-metric.

\begin{rmk}Notice that an asymmetric norm can be defined by the same properties if $X$ is a cone  contained in some vector space $W$. This is needed when defining  the norm of a continuous linear operator between two asymmetric normed spaces (see below).
\end{rmk}

The following example is important in all considerations concerning asymmetric normed spaces.
\begin{example}\label{ex.R-u}
On the field $\mathbb{R}$ of real numbers  the mapping $u(t)=t^+=t\vee 0, \, t\in\mathbb{R},$  is an asymmetric norm with conjugate $\bar u(t)=t^-=(-t)\vee 0$ and the associate norm $u^s(t)=|t|$ is the absolute value of $t$.
\end{example}

\textbf{Linear operators  between asymmetric normed spaces.}

Let $(X,p)$ and $(Y,q)$ be asymmetric normed spaces. A linear operator $A:X\to Y$ is continuous with respect to the topologies $\tau_p$ and $\tau_q$ (we call it $(p,q)$-continuous) if and only if there exists a constant $\beta\ge 0$ such that
\begin{equation}\label{eq1.lin-op-Lip}
q(A(x))\le \beta\, p(x),\end{equation}
for all $x\in X$. Denote by   $\mathcal L_{p,q}(X,Y)$ the set of all continuous linear operators from $X$ to $Y$. Then $\mathcal L_{p,q}(X,Y)$ is a cone in the space $(X,Y)^\#$ of all linear operators from $X$ to $Y$. Actually $\mathcal L_{p,q}(X,Y)$ is contained in the linear space $\mathcal L_s(X,Y)$ of all continuous linear operators from $(X,p^s)$ to $(Y,q^s)$. The number $\beta$ is called a semi-Lipschitz constant for $A$. The smallest of these constants, denoted by $\|A|_{p,q}\,$, is an asymmetric norm
on the cone $ \mathcal L_{p,q}(X,Y)$. It can be also calculated by the formulae:
\begin{equation}\label{eq2.lin-op-norm}
\begin{aligned}
  \|A|_{p,q}&=\sup\{q(A(x)): x\in X, \, p(x)\le 1\}\\
  &=\sup\left\{\frac{q(A(x))}{p(x)}:x\in X,\, p(x)>0\right\}.
  \end{aligned}\end{equation}
  
  In some cases, when the involved asymmetric norms are clear  from the context, we shall omit the subscripts. 

  \subsection{Injective  metric and  quasi-metric spaces}
  Let $(X,\rho)$ be a metric space and $B_\rho[x,\alpha], \, B_\rho[y,\beta]$ two closed balls in $X$. If   $B_\rho[x,\alpha]\cap B_\rho[y,\beta]\ne\emptyset$, then
  \begin{equation}\label{eq1a.hypercv}
  B_\rho[x,\alpha]\cap B_\rho[y,\beta]\ne\emptyset\;\Rightarrow\; \rho(x,y)\le \alpha+\beta.\end{equation}

  If the converse implication also holds, i.e.,
  \begin{equation}\label{eq1b.hypercv}
 \rho(x,y)\le \alpha+\beta\;\Rightarrow\; B_\rho[x,\alpha]\cap B_\rho[y,\beta]\ne\emptyset,
  \end{equation}
  for all $x,y\in X$ and $\alpha,\beta>0,$ then   the metric space $(X,\rho)$ is called \emph{metrically convex}. This is equivalent to  each of the following implications:
  \begin{equation}\label{eq2.hypercv}\begin{aligned}
{\rm (a)}\;\;\rho(x,y)&\le \alpha+\beta \; \Rightarrow\; \exists z\in X,\; d(x,z)\le\alpha\;\mbox{and}\;d(z,y)\le\beta;\\
{\rm (b)}\;\;\rho(x,y)&= \alpha+\beta \; \Rightarrow\; \exists z\in X,\; d(x,z)=\alpha\;\mbox{and}\;d(z,y)=\beta,
  \end{aligned}\end{equation}

  It is easy to check that a normed space is always metrically convex,  because
  \begin{equation*}
  B[x,\alpha]\cap B[y,\beta]\ne\emptyset\iff \|y-x\|\le \alpha+\beta.\end{equation*}

  A metric space  $(X,\rho)$ is called \emph{hyperconvex}  if for any family $x_i,\,i\in I,$ of elements in $X$ and the numbers $r_i>0,\, i\in I,$
  \begin{equation}\label{def.hypercv}
  d(x_i,r_j)\le r_i+r_j\;\mbox{for all}\; i,j\in I\;\Rightarrow\; \bigcap_{i\in I}B[x_i,r_i]\ne \emptyset.\end{equation}

 One says that a metric space $(X,d)$ has the \emph{Binary Intersection Property} (BIP) if any collection of mutually intersecting closed balls has nonempty intersection.   \begin{rmk}\label{re.BIP-compl}
   A  metric space with BIP is complete, a reason for which some authors call these spaces as hypercomplete.
 \end{rmk}

 One shows that a metric space is hyperconvex if and only if it is an absolute nonexpansive retract, i.e., is a nonexpansive retract of every metric space containing it isometrically. A nonexpansive retract of a metric space $(X,d)$ onto a subspace $Y$ is a mapping $R:X\to Y$ such that $R(y)=y$ for all $y\in Y$ and
 $d(R(x),R(x'))\le d(x,x')$ for all $x,x'\in Y.$

 A metric space $(Y,d)$ is called \emph{injective} if for any metric space $(X,\rho)$  and any subspace $Z$ of $X$, every nonexpansive mapping $f:Z\to Y$ admits a nonexpansive extension $F:X\to Y.$  As it was proved in \cite{aronsz-pan56} a metric space is injective if and only if it is hyperconvex.

 Isbell \cite{isbel64a,isbel64b,isbel69} considered injective hulls for metric spaces. The \emph{injective hull} of a metric space $X$ is an injective metric space $Y$ containing $X$ isometrically and isometrically contained  in any injective metric space containing $X$ isometrically (in other words, it is the minimal injective metric space containing $X$ isomatrically).

 Hyperconvex metric spaces were introduced by Aronszajn and Panitchpakdi  \cite{aronsz-pan56} and their study was continued by many mathematicians, in connection with fixed point theory, the extension of mappings and other applications. Good surveys are given in the papers by  Esp\'{\i}nola and Khamsi \cite{Espin-Kham},
 Esp{\'{\i}}nola and Fern{\'a}ndez-Le{\'o}n \cite{Espin-Leon} and in the books by Borkowski \cite{Borkowski}  and Khamsi and Kirk \cite[Chapter 4]{KK}.

 Hyperconvexity in quasi-metric spaces was introduced by Salbany \cite{salbany89}. In our presentation we follow the papers \cite{otaf-ku-con17} and \cite{otaf-ku-kem12}.

A  quasi-metric space  $(Y,d)$ is called:
\begin{itemize}
\item \emph{Isbell-convex} if for any family $y_i,\, i\in I,$ of elements in $Y$ and the numbers $r_i,s_i>0,\, i\in I,$ the following implication holds:
\begin{equation}\label{def.Isb-cv}
  d(y_i,y_j)\le r_i+s_j\;\mbox{ for all }\; i,j\in I\;\;\Rightarrow\;\; \bigcap_{i\in I}\big(B_d[x_i,r_i]\cap B_{\bar d}[y_j,s_j]\big)\ne \emptyset;\end{equation}
\item \emph{metrically convex} if
\begin{equation}\label{def.asym-cv}
d(x,y)\le r+s\; \Rightarrow\; \exists z\in Y,\; d(x,y)\le r\;\;\mbox{and}\;\; d(z,y)\le s, \end{equation}
for all  $x,y\in X$ and $r,s>0$;
\item \emph{Isbell-complete}  if for any family $y_i,\, i\in I,$ of elements in $Y$ and the numbers $r_i>0,\, i\in I,$ the following implication holds:
\begin{equation}\label{def.Isb-co}
  B_d[x_i,r_i]\cap B_{\bar d}[y_j,s_j]\ne \emptyset\; \mbox{ for all }\; i,j\in I\;\;\Rightarrow\;\; \bigcap_{i\in I}\big(B_d[x_i,r_i]\cap B_{\bar d}[y_j,s_j]\big)\ne \emptyset.
  \end{equation}\end{itemize}

  In \cite{otaf-ku-kem12} Isbell-completeness was called \emph{mixed Binary Intersection Property}. The term Isbell-complete was proposed in \cite{otaf-ku-con17}.
  We mention the following properties of  quasi-metric space.
  \begin{prop}[\cite{otaf-ku-kem12}]\label{p1.Isb-cv} Let $(Y,d)$ be a quasi-metric space. Then
 \begin{enumerate}\item  If  $(Y,d) $ is Isbell-convex, then it is Isbell-complete.
 \item The space  $(Y,d)$  is Isbell-convex if and only if  it is Isbell-complete and  metrically convex.
 \item If  $(Y,d) $ is Isbell-convex (Isbell-complete, metrically convex), then the space $(Y,\bar d)$ has the same properties.
 \item If the quasi-metric space $(Y,d) $ is Isbell-convex or Isbell-complete  then the metric space $(Y,d^s)$ is hyperconvex or has the binary intersection property, respectively.
 This implication does not hold for the metric convexity.
 \end{enumerate}\end{prop}

The following property follows from assertion 4 of Proposition \ref{p1.Isb-cv} and  Remark \ref{re.BIP-compl}.
\begin{corol}\label{c.Isb-bicompl}
An Isbell-complete  quasi-metric space  is  bicomplete.
\end{corol}

Recall that a quasi-metric space $(Y,d)$ is called \emph{bicomplete} if the associated metric space $(Y,d^s)$ is complete.

 In the case of an asymmetric normed space $(Y,q)$  the following equivalences hold
\begin{equation}\label{eq.Isb-asym-ns}\begin{aligned}
 q(y_2-y_1)\le r_1+r_2&\iff B_q[y_1,{r_1}]\cap B_{\bar{q}}[y_2,{r_2}]\ne\emptyset\\
 &\iff B_q[y_1,{r_2}]\cap B_{\bar{q}}[y_2,{r_1}]\ne\emptyset,
\end{aligned}
\end{equation}
for any $y_1,y_2\in Y$ and $r_1,r_2>0.$

  Indeed, if $y\in  B_q[y_1,{r_1}]\cap B_{\bar{q}}[y_2,{r_2}]$, then
  $$
  q(y_2-y_1)\le q(y_2-y)+q(y-y_1)\le r_2+r_1\,.$$

  A similar conclusion can be drawn taking an $y$ in $ B_q[y_1,{r_2}]\cap B_{\bar{q}}[y_2,{r_1}].$

  Suppose now that $q(y_2-y_1)\le r_1+r_2.$  Then
  $$
  y:=\frac{r_2}{r_1+r_2}y_1+\frac{r_1}{r_1+r_2}y_2\in B_q[y_1,{r_1}]\cap B_{\bar{q}}[y_2,{r_2}],$$
  and
  $$
 y':=\frac{r_1}{r_1+r_2}y_1+\frac{r_2}{r_1+r_2}y_2\in B_{q}[y_1,{r_2}]\cap B_{\bar q}[y_2,{r_1}].$$

 The equivalences \eqref{eq.Isb-asym-ns} imply:
 \begin{corol}\label{c.Isb-cv-co-asym-ns}
 An asymmetric normed space $(Y,q)$ is Isbell-convex if and only if it is Isbell-complete.\end{corol}

 In what follows  we shall use the term \emph{mixed Binary Intersection Property} (mixed BIP) to design Isbell-completeness.

  \section{The extension of linear operators between asymmetric normed spaces}

 The injectivity of an asymmetric normed space $(Y,q)$ is defined as in the case of normed spaces (the existence of norm-preserving extensions of continuous linear operators).

 A typical example of an injective asymmetric normed space is the space $(\ell^\infty(I),\tilde q)$.  For a  nonempty set  $I$ denote by $\ell^\infty(I)$ the linear space of all bounded
  functions $x:I\to\mathbb{R}.$. For $x=(x_i)_{i\in I}\in\ell^\infty(I)$  put
  \begin{equation}\label{def.asym-ell-infty}
  \tilde q(x)=\sup_{i\in I}x_i^+=  \sup_{i\in I}x_i
  \end{equation}
  where $\alpha^+=\max\{\alpha,0\}$ for every $\alpha\in \mathbb{R}.$

  Then $\tilde q$ is an asymmetric norm on $\ell^\infty(I)$.

  In the following proposition we collect some properties of injective asymmetric normed spaces, similar to some well-known properties  of injective normed spaces.
  \begin{prop}\label{p1.ell-infty-inj} Then the following properties hold.
  \begin{enumerate}
  \item Let $(X,p)$ be an asymmetric normed space.  A linear operator $T:X\to \ell^\infty(I)$ is continuous if and only if there exist the  continuous linear functionals $\varphi_i:(X,p)\to (\mathbb{R},u),\, i\in I,$
  such that
  \begin{equation*}
  T(x)_i=\varphi_i(x),\, i\in I,\, x\in X\;\mbox{ and }\; \|T|=\sup_{i\in I}\|\varphi_i|.\end{equation*}
  \item  The asymmetric normed space $(\ell^\infty(I),\tilde q)$ is injective.
  \item  Every asymmetric normed space   can be linearly and isometrically embedded in $(\ell^\infty(I),q)$ for some set $I$, where $q$ is given by \eqref{def.asym-ell-infty}.
      \item Any asymmetric normed space isometrically isomorphic to an injective  asymmetric normed space is  injective  itself.
  \item If  $(Y,q)$  is an injective asymmetric normed space, then  for every asymmetric normed space $X$ containing it  isometrically, there exists a norm-one linear projection $P:X\to Y.$
      \item Conversely, if an  asymmetric normed space    $(Y,q)$  is such that  for every asymmetric normed space $(X,p)$ containing it  isometrically, there exists a norm-one linear projection $P:X\to Y,$  then $(Y,q)$ is injective.
      \item An asymmetric normed space $(Y,q)$ is injective if and only if it is linearly isometric to a complemented subspace of $(\ell^\infty(I)$, for some set $I$.
  \end{enumerate}\end{prop}\begin{pf}
    1. Suppose that $A:X\to  \ell^\infty(I)$ is linear and $(p,\tilde q)$-continuous. Then, for every $i\in I$, the functional $\varphi_i:X\to\mathbb{R}$  given by
    $$
    \varphi_i(x)=T(x)_i,\, x\in X,$$
    is linear and continuous and $\|\varphi_i|\le\|T|.$

    Also
\begin{align*}
    \|T|&=\sup\{\tilde q(A(x)):x\in X,\, p(x)\le 1\}\\
    &=\sup\{\sup_{i\in I}A(x)_i:x\in X,\, p(x)\le 1\}\\&=\sup_{i\in I}\sup\{\varphi_i(x):x\in X,\, p(x)\le 1\}\\&=\sup_{i\in I}\|\varphi_i|.
\end{align*}

    Conversely, if $\varphi_i:X\to \mathbb{R},\, i\in I,$ are given $(p,u)$-continuous linear functionals, then the operator $T:X\to  \ell^\infty(I)$ defined by $T(x)_i=\varphi_i(x),\, i\in I,\, x\in X$, is linear, $(p,\tilde q)$-continuous and
    $\|T|=\sup_{i\in I}\|\varphi_i|.$

2. Let $(X,p)$ be an asymmetric normed space, $Z$ a subspace of $X$ and $T:Z\to\ell^\infty(I)$ a  $(p,\tilde q)$-continuous linear operator. Consider the $(p,u)$-continuous linear  functionals $\varphi_i:Z\to\mathbb{R},\, i\in I,$ representing the operator $T$, as in  (1).  Then, by the Hahn-Banach  extension theorem in asymmetric normed spaces (\cite[Theorem 2.2.2]{Cobzas}, the  functionals $\varphi_i$ admit norm-preserving linear extensions $\psi_i:X\to\mathbb{R}.$ Defining $\tilde T:X\to \ell^\infty(I)$ by $\tilde T(x)_i=\psi_i(x),\, i\in I,\, x\in X,$
it follows that $\tilde T$ is a continuous linear extension of $T$ with norm
$$
\|\tilde T|=\sup_{i\in I}\|\psi_i|= \sup_{i\in I}\|\varphi_i|=\|T|.$$

3. Let $\gamma$ be the density characteristic of $(X,p^s)$, i.e., the least cardinal number of a dense subset of    $(X,p^s)$, and let $\{x_i:i\in I\}$  be a dense subset of $(X,p^s)$ with card\,$(I)=\gamma$.

Appealing again to \cite[Theorem 2.2.2]{Cobzas}, for every $i\in I$ with $p(x_i)>0$ there exists a continuous linear functional $\varphi_i:(X,p)\to(\mathbb{R},u)$ with
$$
\|\varphi_i|=1\;\mbox{ and }\; \varphi(x_i)=p(x_i).$$

If $p(x_i)=0$, then take $\varphi_i\equiv 0.$

Define the operator $A:X\to  \ell^\infty(I)$ by
$$
A(x)_i=\varphi_i(x),\, i\in I,\, x\in X.$$

Then
$$
\tilde q(A(x))=\sup_{i\in I} A(x)_i= \sup_{i\in I} \varphi_i(x)\le p(x),$$
for all $x\in X$.  This inequality shows that the operator $A$ is $(p,\tilde q)$-continuous.

For $x\in X$, let $(x_{i_k})_{k\in\mathbb{N}}$ be a sequence $p^s$-convergent to $x$, that is,
$$
x_{i_k}\xrightarrow{p}x\;\mbox{ and }\;  x_{i_k}\xrightarrow{\bar p}x,$$
as $k\to\infty.$

We have
$$
\tilde q(A(x_{i_k}))\ge A(x_{i_k})_{i_k}=\varphi_{i_k}(x_{i_k})=p(x_{i_k}),$$
for all $k\in \mathbb{N}.$  By Proposition \ref{p.top-qsm1}) $\tilde q$ is $\tau(\tilde q)$-usc   and $p$ is $\tau(\bar p)$-lsc, so that, taking into account the $(p,\tilde q)$-continuity of $A$,
 one obtains
$$
\tilde q(A(x))\ge\limsup_{k\to\infty}\tilde q(A(x_{i_k})\ge\liminf_{k\to\infty}p(x_{i_k})\ge p(x).$$

It follows that $\tilde q(A(x))=p(x)$ for all $x\in X.$

4. Let $(Y,q)$ and  $(\tilde Y,\tilde q)$ be asymmetric normed spaces, with $Y$ injective, and let $J:Y\to \tilde Y$ be an isometric isomorphism. If $Z$ is a subspace of an asymmetric normed space $(X,p)$ and $T:Z\to \tilde Y$ is  a $(p,\tilde q)$-continuous linear operator, then  $J^{-1}\circ T:Z\to Y$ is linear and $(p,q)$-continuous, so that, by the injectivity of  $(Y,q)$, it admits a norm-preserving linear extension $S:X\to Y$. But then $J\circ S:X\to \tilde Y$ is a norm preserving linear  extension of $T$, showing that $(\tilde Y,\tilde q)$ is injective too.

5.  Suppose that $(Y,q)$ is an injective asymmetric normed space contained in an asymmetric normed space $(X,p)$ with $p|_Y=q.$ Then the identity operator 
${\rm Id}_Y:(Y,p)\to (Y,q)$  is linear, $(p,q)$-continuous  and of  norm one. By the injectivity of $Y$ it admits a norm-preserving linear extension $P:(X,p)\to (Y,q)$ which is a norm-one linear projection  of $X$ onto $Y$.

6.  Suppose that the asymmetric normed space  $(Y,q)$ satisfies the hypothesis of (5). Let $Z$ be a subspace of an asymmetric normed space  $(X,p)$ and $T:Z\to Y$ a $(p,q)$-continuous linear operator.
By  assertion 3, $\,(Y,q)$ is  isometrically isomorphic to a subspace of $(\ell^\infty(I),\tilde q)$ for some set $I$, where $\tilde q$ is given by \eqref{def.asym-ell-infty}, i.e.,
there exists a linear isometry $J:Y\to\ell^\infty(I).$  Since the asymmetric normed space $\ell^\infty(I)$ is injective, the operator $J\circ T$ admits a norm preserving extension $S:X\to \ell^\infty(I).$  By 4 $J(Y)$ is also injective, so that, by 5,  there exists a norm-one linear projection $P:\ell^\infty(I)\to J(Y).$ Let $j:J(Y)\to Y$ be the inverse of the mapping $J:Y\to J(Y),$
that is $j(J(y))=y$ for all $y\in Y$. Then $j$ is an isometric isomorphism and $\tilde T=j\circ P\circ S$ is a norm-preserving linear extension of $T$.

7. This follows from the assertions 4 and 5 of the proposition, taking into account that a subspace $Y$  of an asymmetric normed space $X$ is complemented   in $X$ if and only if there exists a continuous linear projection $P:
X\to Y.$
 \hfill $\Box$\end{pf}

 \begin{thm}
 \label{t1.asym-Nachbin}
An asymmetric normed space  $(Y,q) $  has  the mixed  BIP if and only if for any asymmetric normed space $(X,p)$ and  any $p$- and $\bar p$-closed subspace $Z$ of $X$
 every $(p,q)$-continuous linear operator $T:Z\to Y$ admits a norm preserving extension  to $X$, that is there exists  a $(p,q)$-continuous linear operator $S:X\to Y$
such that $S|_Z=T$ and $\|S|_{p,q}=\|T|_{p,q}.$
  \end{thm}
\begin{pf}  \emph{Sufficiency}.\,  Suppose that $(Y,q)$ is an asymmetric normed space satisfying the mixed BIP and  let   $(X,p)$ be another asymmetric normed space, $Z$ a subspace of $X$ that is both $p$- and $\bar p$-closed and let $T:Z\to Y$  be a $(p,q)$-continuous linear operator. Without restricting the generality we can suppose
$\|T|_{p,q}=1,$ i.e.,
$$
q(T(z))\le p(z)$$
for all $z\in Z$.

Let $x_0\in X\smallsetminus Z$ and  $U=T(Z)\subset Y$.  For $u\in U$ put
\begin{equation}\label{def.r(u)}\begin{aligned}
  &r(u)=d_p(x_0,T^{-1}(u))=\inf\{p(z-x_0) : z\in Z,\, T(z)=u\}, \mbox{ and}\\
   &s(u)=d_{\bar p}(x_0,T^{-1}(u))=\inf\{p(x_0-z) : z\in Z,\, T(z)=u\}
\end{aligned}\end{equation}

Notice that $r(u)>0$ and $s(u)>0$, for all $u\in U.$ Indeed, if for some $u\in U,\, r(u)=0,$ then there exists a sequence $(z_n)$ in $Z$ such that $T(z)_n=u,\, n\in\mathbb{N},$ and
$p(z_n-x_0)\to 0$, implying that $x_0\in \clos_p(Z)=Z$, a contradiction.  Similarly, $s(u)=0$ would imply $x_0\in \clos_{\bar p}(Z)=Z,$ a contradiction again.

We have\smallskip

\textbf{Claim 1. } ${\ds q(u_2-u_1)\le r(u_2)+s(u_1)}$ for all $u_1,u_2\in U$.\smallskip

Indeed,
$$
q(u_2-u_1)=q(T(z_2-z_1))\le p(z_2-x_1)\le p(z_2-x_0)+p(x_0-z_1)\,,$$
for all $z_1,z_2\in Z$ such that $T(z)_1=u_1$ and  $T(z)_2=u_2.$ Passing to infimum with respect to these $z_1,z_2,$ one obtains
$$
q(u_2-u_1)\le d_p(x_0,T^{-1}(u_2))+  d_{\bar{p}}(x_0,T^{-1}(u_1))=r(u_2)+s(u_1)\,.$$

Now, by   Claim 1 and  \eqref{eq.Isb-asym-ns},
$$
B_q[u_1,{s(u_1)}]\cap  B_{\bar{q}}[u_2,{r(u_2)}]\ne\emptyset\,,$$
for all $u_1,u_2\in U$, so that, by hypothesis, there exists $y_0\in Y$ such that
\begin{equation}\label{def.y}
y_0\in\bigcap_{u\in U}\left(B_q[u,{s(u)}]\cap  B_{\bar{q}}[u,{r(u)}]\right)\,.\end{equation}

Define now the operator $S:Z\dotplus\mathbb{R} x_0\to Y $ by
$$
S(z+t x_0)=T(z)+t y_0$$
for all $z\in Z$ and $t\in\mathbb{R}$.

The operator $S$ is obviously linear, $S(x_0)=y_0$   and $S|_Z=T$, so it remains to show that it is continuous and\, $\|S|_{p,q}=1$, which is equivalent to
\begin{equation}\label{ineq.S1}
q(S(z+t x_0))=q(T(z)+t y_0)\le p(z+t x_0)\,,
\end{equation}
for all $z\in Z$ and $t\in\mathbb{R}$.

Let us show first that
\begin{equation}\label{ineq.S2}\begin{aligned}
&{\rm (a)}\; \quad q(T(z)+y_0)\le p(z+x_0)\;\mbox{ and}\\
&{\rm (b)}\;\quad  q(T(z)-y_0)\le p(z-x_0),
\end{aligned}\end{equation}
for all $z\in Z$.

Since $y_0\in B_{\bar{q}}[T(z),{r(T(z))}]$ it follows
$$
q(T(z)-y_0)\le r(T(z))\le p(z-x_0)\,,$$
that is, \eqref{ineq.S2}.(b) holds.

To prove (a), observe that $y_0\in B_{q}[T(-z),{s(T(-z))}]$, so that
\begin{align*}
  q(T(z)+y_0)&=q(y_0-T(-z))\le s(T(-z))\\&\le p(x_0-(-z))=p(x_0+z)\,,
\end{align*}
showing that  \eqref{ineq.S2}.(a) holds too.

 Now, if $t >0$, then, by  \eqref{ineq.S2}.(a),
  $$
  q(S(z+t x_0))=q(Tx+t y_0)=t q(T(t^{-1}z)+y_0)\le t p(t^{-1}z+x_0)=p(z+t x_0)\,.$$

  If $t <0$, then, by  \eqref{ineq.S2}.(b),
  $$
  q(S(z+t x_0))=q(Tx+t y_0)=-t q(T(-t^{-1}z)-y_0)\le -t p(-t^{-1}z-x_0)=p(z+t x_0)\,.$$

  Consequently, $S$ satisfies \eqref{ineq.S1}.
  \hfill $\Box$\end{pf}

  \begin{rmk}
    A subset $Z$ of an asymmetric normed space $(X,p)$ that is  $p$- or $\bar p$-closed is also $p^s$-closed, but a $p^s$-closed subset of $X$ need not be
    $p$-closed nor $\bar p$-closed.
  \end{rmk}

Indeed, suppose that $Z$ is $p$-closed and let  $(z_n)$ be a sequence in $Z$,  $p^s$-convergent to some $x\in X.$  Then $z_n\xrightarrow{p}x,$   so that $x\in Z$, showing that $Z$ is $p^s$-closed.
The case of the conjugate asymmetric norm $\bar p$ can be treated similarly. 

The converse implication could fail. For instance, in $(\mathbb{R},u)$ (see Example \ref{ex.R-u}), the ball
$B_{u^s}[0,1]=[-1,1]$ is $u^s$-closed  but  neither $u$-closed nor $\bar u$-closed.   Indeed, for every $r>0,$ $B_{u^s}[0,1]=[-1,1]\subset(-\infty,2+r)=B_u(2,r)0, $ so that there is no a $u$-neighborhood  of 2, disjoint from $[-1,1]$. Hence $B_{u^s}[0,1]$ is  not
$u$-closed.  Similarly, for every $r>0,\, B_{\bar u}(-2,r)=(-2-r,\infty)\supset [-1,1]=B_{u^s}[0,1]$, showing that $B_{u^s}[0,1]$ is not $\bar u$-closed.

 The proof of the necessity part of Theorem \ref{t1.asym-Nachbin} is  more elaborate  and requires some preliminaries.

 Let $(Y,p)$ be an asymmetric normed space and $Z$ a   proper subspace of $Y$. For $\xi\in Y\smallsetminus Z$ consider the functions $\rho_1,\rho_2:Z\to\mathbb{R}_+,$ given by
\begin{equation}\label{def.rho}
   \rho_1(z)=p(z-\xi)\;\mbox{and}\; \rho_2(z)=p(\xi-z)\;\mbox{for all}\; z\in Z.
\end{equation}

 It is easy to check that the functions $\rho_1,\rho_2$ have the following properties:
 \begin{itemize}
 \item[(a)] $\rho_1(z')+\rho_2(z)\ge p(z'-z)$;
 \item[(b1)] $\rho_1(z')-\rho_2(z)\le p(z'-z)$;
 \item[(b2)] $\rho_2(z')-\rho_1(z)\le p(z-z')$;
 \item[(c)] the functions $\rho_1,\rho_2$ are convex;
 \item[(d)] $\rho_1(z)+\rho_2(z)>0,$
   \end{itemize}
   for all $z,z'\in Z.$

   We shall show that any pair  of functions  $\rho_1,\rho_2:Z\to \mathbb{R}_+$  satisfying (a)-(d) are of the form \eqref{def.rho}, for an appropriate asymmetric norm.
  \begin{prop}\label{p1.asym-rho}
   Let $(Z,p)$ be an asymmetric normed space.  \begin{enumerate}\item
     Suppose that $Z$ is contained in a vector space $W$,  $\xi\in W\smallsetminus Z,$  and $Y=Z\dotplus\mathbb{R} \xi$  (direct algebraic sum).
 If $\rho_1,\rho_2:Z\to \mathbb{R}_+$ are two functions on $Z$  satisfying the conditions (a)-(d), then the mapping $\mu:Y\to\mathbb{R}_+$ defined
   for $y=z+t\xi\in Y$ by
\begin{equation}\label{def.asym-mu}
   \mu(y)=\begin{cases}
     p(z)\quad &\mbox{if}\quad t=0\\
     t\rho_2\left(\frac{-z}t\right)  \quad &\mbox{if}\quad t>0,\\
     -t\rho_1\left(\frac{-z}t\right)  \quad &\mbox{if}\quad t<0,
   \end{cases}\end{equation}
      is an asymmetric norm  on $Y$ and
$$
\rho_1(z)=\mu(z-\xi),\;\; \rho_2(z)=\mu(\xi-z),$$
for all $z\in Z.$
\item If the functions  $r_1,r_2: Z\to\mathbb{R}_+$ satisfy the condition (a), then there exist two functions $\rho_1,\rho_2:Z\to\mathbb{R}_+$ satisfying the conditions (a)-(c) such that  $r_1\ge \rho_1,\, r_2\ge \rho_2.$
    \end{enumerate}\end{prop}    \begin{pf}
1.  Observe that, by the definition of $\mu$,
   $$
   \mu(z-\xi)=\rho_1(z)\;\mbox{ and }\; \mu(\xi-z)=\rho_2(z),$$
   for all $z\in Z.$

   Let $y=z+t\xi\ne 0$. Then $-y=-z+(-t)\xi.$ If $t=0$, then $z\ne 0$ and
   $$\mu(y)+\mu(-y)=p(z)+p(-z)>0.$$

   If $t\ne 0$, then, by (d),
   $$
   \mu(y)+\mu(-y)=t\rho_2\left(\frac{-z}t\right) +t\rho_1\left(\frac{-z}t\right)>0,$$
   for $t>0$, and
    $$
   \mu(y)+\mu(-y)=(-t)\rho_1\left(\frac{-z}t\right) +(-t)\rho_1\left(\frac{-z}t\right)>0,$$
   for $t<0$.\smallskip

\emph{   The positive homogeneity.}

Let $\alpha>0$ and $t>0$. Then
$$
\mu(\alpha y)=\alpha t\rho_2\left(\frac{-\alpha z}{\alpha t}\right)=\alpha\mu(y).$$

A similar calculation can be done if $t<0$ (working with $\rho_1$).\smallskip

\emph{The triangle  inequality.}

  Let
$$y=z+t\xi,\; y'=z'+t'\xi,\; y+y'=z+z'+(t+t')\xi$$
be elements in $Y$.

We have to prove that
\begin{equation}\label{ineq.asym-tr}
\mu(y+y')\le  \mu(y)+\mu(y').\end{equation}

To this end  we have to consider several cases.
\smallskip

I.\; $t=t'=0.$

In this case, \eqref{ineq.asym-tr} is equivalent to
$$
p(z+z')\le p(z)+p(z'),$$
which is true since $p$ is an asymmetric norm.\smallskip

II. \; \emph{Only one of the numbers} $t,t'$ \emph{is null, say} $t'=0.$

Then $y+y'=z+z'+t\xi,$ so that \eqref{ineq.asym-tr} is equivalent to
$$
t\left(\rho_2\left(\frac{-z-z'}{t}\right)-\rho_2\left(\frac{-z}{t}\right)\right)\le p(z'),$$
if $t>0$, which is true by (b2), and to

$$
(-t)\left(\rho_1\left(\frac{-z-z'}{t}\right)-\rho_1\left(\frac{-z}{t}\right)\right)\le p(z'),$$
if $t<0,$ which follows from (b1).\smallskip

III.\; \emph{Both numbers $t,t'$ are nonnull}.

We split it into several subcases.\smallskip

III.(i) \; $t,t'>0$.

In this case \eqref{ineq.asym-tr} is equivalent to
$$
\rho_2\left(\frac{-z-z'}{t+t'}\right)\le \frac{t}{t+t'}\cdot\rho_2\left(\frac{-z}{t}\right)+\frac{t'}{t+t'}\cdot\rho_2\left(\frac{-z'}{t'}\right),$$
which follows from the convexity of $\rho_2$ (condition (c)) and the equality
\begin{equation}\label{eq1.asym-tr}
\frac{-z-z'}{t+t'}= \frac{t}{t+t'}\cdot\frac{-z}{t}+\frac{t'}{t+t'}\cdot\frac{-z'}{t'}.\end{equation}

III.(ii)\; $t,t'<0$.

In this case \eqref{ineq.asym-tr} is equivalent to
$$
\rho_1\left(\frac{z+z'}{-(t+t')}\right)\le \frac{-t}{-(t+t')}\cdot\rho_1\left(\frac{z}{-t}\right)+\frac{-t'}{-(t+t')}\cdot\rho_1\left(\frac{z'}{-t'}\right),$$
which follows from the convexity of $\rho_1$ and the equality \eqref{eq1.asym-tr}.\smallskip

III.(iii)\; $t>0,\, t'=-t.$

In this case $y+y'=z+z'$, so that \eqref{ineq.asym-tr} is equivalent to
$$
p(z+z')\le t\left(\rho_2\left(\frac{-z}{t}\right)+\rho_1\left(\frac{z'}{t}\right)\right),$$
which follows from (a).\smallskip

III.(iv)\; $t>0,\, t'<0,\, t+t'>0$.

In this case \eqref{ineq.asym-tr} is equivalent to
$$
(t+t')\rho_2\left(\frac{-z-z'}{t+t'}\right)\le t\rho_2\left(\frac{-z}{t}\right)-t'\rho_1\left(\frac{-z'}{t'}\right),$$
which, in its turn, is  equivalent to
\begin{align*}
  A&:=(t+t')\left(\rho_2\left(\frac{-z-z'}{t+t'}\right)-\rho_2\left(\frac{-z}{t}\right)\right)\\
  &\le (-t')\left(\rho_2\left(\frac{-z}{t}\right)+\rho_1\left(\frac{-z'}{t'}\right)\right)=:B.
\end{align*}

But, by (b2) and (a),
$$
A\le p\left(z'-\frac{t'}{t}z\right)\le B.$$

III.(v)\; $t>0,\, t'<0,\, t+t'<0.$

By a similar  calculation  to that done  in the case III.(iv), one shows that \eqref{ineq.asym-tr} is  equivalent to
\begin{align*}
  A'&:=-(t+t')\left(\rho_1\left(\frac{z+z'}{-(t+t')}\right)-\rho_1\left(\frac{-z'}{t'}\right)\right)\\
  &\le t\left(\rho_2\left(\frac{-z}{t}\right)+\rho_1\left(\frac{-z'}{t'}\right)\right)=:B'.
\end{align*}

By (b1) and (a),
$$
A'\le p\left(z-\frac{t}{t'}z'\right)\le B'.$$

In conclusion, the triangle inequality \eqref{ineq.asym-tr} is proved in all possible cases, concluding the proof that $\mu$ is an asymmetric norm on $Y$.\smallskip

2. Consider the set $\Lambda$ of all pairs $f=(f_1,f_2) $ of functions $f_1,f_2:Z\to\mathbb{R}_+$ satisfying (a).  The set $\Lambda$ is nonempty - for instance,
$f=(f_1,f_2)\in\Lambda$,  where $f_1(z)=f_2(z)=p(z), \, z\in Z.$

Consider on $\Lambda$ the pointwise order given, for $f=(f_1,f_2)$ and $g=(g_1,g_2)$ in $\Lambda$, by
$$
f\le g\iff f_k(z)\le g_k(z),\, k=1,2,$$
for all $z\in Z$.

Then $\Lambda$ is inductively ordered with respect to this order, i.e., every totally ordered subset of $\Lambda$ has a lower bound. Indeed, if $\Gamma\subseteq \Lambda$ is totally ordered, then   $\varphi=(\varphi_1,\varphi_2)$  is a lower bound for $\Gamma,$  where
$$
\varphi_k(z)=\inf\{f_k(z):f\in\Gamma\},\,z\in Z, \; k=1,2,$$
We have to check only  that $\varphi$ belongs to $\Lambda$. Let $z,z'\in Z$. Then, for every $\varepsilon>0$ there exist
 $f=(f_1,f_2)$ and $g=(g_1,g_2)$ in $\Gamma$ such that
$$
  f_1(z)\le \varphi_1(z)+\varepsilon\quad\mbox{and}\quad   g_2(z')\le\varphi_2(z')+\varepsilon.
$$

Since $\Gamma$ is totally ordered, the pairs  $f,g$ are comparable, say $f\le g.$   Then
$$
p(z-z')\le f_1(z) +f_2(z')\le\varphi(z)+\varphi(z')+2\varepsilon.$$

Since $\varepsilon>0$ was arbitrarily chosen, it follows $p(z-z')\le \varphi_1(z)+\varphi_2(z'),$ that is $\varphi\in\Lambda.$

By Zorn's Lemma, there exists a minimal element $\rho=(\rho_1,\rho_2)\in \Lambda.$

We show that, for every $z\in Z,$
\begin{equation}\label{eq.BIP-min-asym}\begin{aligned}
\rho_1(z)&=0\vee\sup\{p(z-z')-\rho_2(z'): z'\in Z,\, z'\ne z\},\;\mbox{and}\\
\rho_2(z)&=0\vee\sup\{p(z'-z)-\rho_1(z'): z'\in Z,\, z'\ne z\}.
\end{aligned}
\end{equation}

For $z\in Z$ consider the functions $\varphi_1,\varphi_2:Z\to \mathbb{R}_+$ defined by $\varphi_k(z')=\rho_k(z')$, for  $z'\in Z\smallsetminus\{z\},\, k=1,2,$
and

\begin{equation*}\begin{aligned}
\varphi_1(z)&=0\vee\sup\{p(z-z')-\rho_2(z'): z'\in Z,\, z'\ne z\},\\
\varphi_2(z)&=0\vee\sup\{p(z'-z)-\rho_1(z'): z'\in Z,\, z'\ne z\}.
\end{aligned}\end{equation*}

Observe that, since $\rho$  satisfies (a),
\begin{align*}
  p(z-z')-\rho_2(z')&\le \rho_1(z),\;\mbox{and}\\
   p(z'-z)-\rho_1(z')&\le  \rho_2(z),
\end{align*}
for all $z'\in Z\smallsetminus\{z\}$. Hence the considered suprema are finite,  $\varphi_1$ and $\varphi_2$ are well-defined  and

\begin{equation}\label{eq.rho-phi}
\varphi_1(z)\le\rho_1(z),\quad \varphi_2(z)\le \rho_2(z).\end{equation}

We have \begin{align*}
  \varphi_1(z)&\ge p(z-z')-\rho_2(z')= p(z-z')-\varphi_2(z'),\;\mbox{and}\\
  \varphi_2(z)&\ge p(z'-z)-\rho_1(z')= p(z'-z)-\varphi_1(z'),
  \end{align*}  so that
  $$
  \varphi_1(z)+\varphi_2(z')\ge p(z-z')\;\mbox{ and }\; \varphi_1(z')+\varphi_2(z)\ge p(z'-z),$$
   for all $z'\ne z.$

   If $z',z''\in Z\smallsetminus \{z\}$,  then
$$
\varphi_1(z')+\varphi_2(z'') =\rho_1(z')+\rho_2(z'')\ge p(z'-z''),$$
because $\rho=(\rho_1,\rho_2) $ belongs to $\Lambda.$

All these show that $\varphi=(\varphi_1,\varphi_2)\in\Lambda.$

Since $\rho=(\rho_1,\rho_2)$ is minimal, it follows $\rho_k(z)\le\varphi_k(z),\, k=1,2,$
so that, taking into account \eqref{eq.rho-phi},  $$\rho_1(z)=\varphi_1(z)\;\mbox{ and }\;  \rho_2(z)=\varphi_2(z).$$

This shows that the equalities \eqref{eq.BIP-min-asym} hold.\smallskip

\emph{Property} (b1).  Since
\begin{align*}
  \rho_1(z)\ge p(z-z'')-\rho_2(z'')\ge p(z'-z'')-\rho_2(z'')-p(z'-z),
\end{align*}
for all $z''\in Z,$ it follows  that
$$\rho_1(z)-\rho_1(z')\ge-p(z'-z)\iff \rho_1(z')-\rho_1(z)\le p(z'-z).$$

\emph{Property} (b2).  Since
\begin{align*}
  \rho_2(z)\ge p(z''-z)-\rho_1(z'')\ge p(z''-z')-\rho_1(z'')-p(z-z'),
\end{align*}
for all $z''\in Z,$ it follows  that
$$\rho_2(z)-\rho_2(z')\ge-p(z-z')\iff \rho_2(z')-\rho_2(z)\le p(z-z').$$

\emph{Property} (c)\, (\emph{the convexity}).

Let $z,z'\in Z$ and $\alpha\in[0,1]. $ We have  prove that
\begin{equation}\label{eq.cv-asym}
\rho_1(\alpha z+(1-\alpha)z')\le \alpha \rho_1(z)+(1-\alpha)\rho_1(z').\end{equation}

For every $z''\in Z,$
\begin{align*}
\alpha \rho_1(z)+(1-\alpha)\rho_1(z')&\ge \alpha \big(p(z-z'')-\rho_2(z'')\big)+(1-\alpha) \big(p(z'-z'')-\rho_2(z'')\big)\\
&\ge p(\alpha(z-z'')+(1-\alpha)(z'-z''))-\rho_2(z'')\\
&=p(\alpha z+(1-\alpha)z'-z'')-\rho_2(z'').
\end{align*}

Passing to supremum with respect to $z''\in Z$, one obtains \eqref{eq.cv-asym}.

The proof of the convexity of $\rho_2$ is similar.\hfill $\Box$\end{pf}

Now we can proceed  to the proof of the necessity part of Theorem \ref{t1.asym-Nachbin}.

\begin{pf}[Proof of Theorem \ref{t1.asym-Nachbin}. Necessity]
  The proof will be given by contradiction.

Suppose that, for an asymmetric normed space $(Z,q)$,  there exist the elements $x_i\in Z$ and the numbers $r_i,s_i>0,\, i\in I,$ such that
\begin{equation*}
B_q[x_i,r_i]\cap B_{\bar q}[x_j,s_j]\ne\emptyset,\end{equation*}
for all $i,j\in I,$ but
\begin{equation}\label{eq.asym-ety}
\bigcap_{i\in I}\left(B_q[x_i,r_i]\cap B_{\bar q}[x_i,s_i]\right)=\emptyset.
\end{equation}

Denote by $A$ the set of all centers of these balls and, for $x\in A$, let
$$
I_x=\{i\in I:x_i=x\}.$$

For $x\in A,$ let
$$
r_1(x)=\inf\{s_i:i\in I_x\}\;\mbox{ and }\;  r_2(x)=\inf\{r_i:i\in I_x\}.$$

Let $x,x'\in A,\, x\ne x'.$ Then
$$
q(x'-x) =q(x_j-x_i)\le s_j+r_i,$$
for all $i\in I_x$ and $j\in I_{x'}.$ Passing to infimum with respect to these $i,j$, one obtains

$$q(x'-x) \le r_1(x')+r_2(x).$$

By  \eqref{eq.Isb-asym-ns}, this implies
$$
B_q[x,r_2(x)]\cap B_{\bar q}[x',r_1(x')]\ne \emptyset.$$

Fix an element $\bar x\in A$ and, for $z\in Z\smallsetminus A$, put
$$
r_1(z)=q(z-\bar x)+r_1(\bar x)\;\mbox{ and }\; r_2(z)=q(\bar x-z)+r_2(\bar x).$$

In this way, the functions $r_1,r_2$ are defined on the whole $Z$.

We shall show that, for all $z,z'\in Z,$
\begin{equation}\label{eq2a.asym-int-balls}
B_q[z,r_2(z)]\cap B_{\bar q}[z',r_1(z')]\ne \emptyset,\end{equation}
which is equivalent to
\begin{equation}\label{eq2.asym-int-balls}
q(z'-z)\le r_1(z')+r_2(z).\end{equation}

Obviously, \eqref{eq2.asym-int-balls} holds for all $z,z'\in A.$

Take $z=x\in A$ and $z'\in Z\smallsetminus A$. Then
$$
q(z'-x)\le q(z'-\bar x)+q(\bar x-x)\le q(z'-\bar x)+r_1(\bar x)+r_2(x)=r_1(z')+r_2(x),$$
so that \eqref{eq2.asym-int-balls}  holds.

One shows, similarly, that if $z\in Z\smallsetminus A $ and $z'=x\in A,$ then
$$q(x-z)\le r_1(x)+r_2(z),$$
so that \eqref{eq2.asym-int-balls} holds in this case too.

If $z,z'\in Z\smallsetminus A$, then
\begin{align*}
  r_1(z')+r_2(z)&=q(z'-\bar x)+r_1(\bar x)+q(\bar x-z)+r_2(\bar x)\\
  &\ge q(z'-z),
\end{align*}
showing that \eqref{eq2.asym-int-balls} also holds  in this case.

By Proposition \ref{p1.asym-rho}, there exists a pair of  functions $\rho=(\rho_1,\rho_2)\le (r_1,r_2)$ satisfying the conditions (a)-(c).

Let us show that $\rho$ also satisfies the condition (d). Indeed, if, for some $z\in Z,$
$$\rho_1(z)=\rho_2(z)=0,$$
then, by (a),   for every $z'\in Z,$
$$
q(z'-z)\le\rho_1(z')+\rho_2(z)\le r_1(z'),$$
showing that $z\in B_{\bar q}[z',r_1(z')]$.

Also,

$$
q(z-z')\le\rho_1(z)+\rho_2(z')\le r_2(z'),$$
showing that $z\in B_{q}[z',r_2(z')]$.

Hence, $$z\in\bigcap_{z'\in Z}\left(B_{q}[z',r_2(z')]\cap B_{\bar q}[z',r_1(z')]\right),$$
in contradiction to \eqref{eq.asym-ety}.

Suppose now that $Z$ is a proper subspace of some vector space $W, \; \xi\in W\smallsetminus Z,$ and $Y=Z\dotplus \mathbb{R}\xi.$
Then, by Proposition \ref{p1.asym-rho}.1, there exists an asymmetric norm $\mu$ on $Y$, given by \eqref{def.asym-mu}, such that
$$
\mu(z-\xi)=\rho_1(z)\;\mbox{ and }\; \mu(\xi-z)=\rho_2(z),$$
for all $z\in Z.$

Suppose that there exists a norm-one linear projection $P:Y\to Z$ and let $z_0=P(\xi).$ Then, for every $z\in Z,$
$$
q(z-z_0)=\mu(z-z_0)=\mu(P(z-\xi))\le\mu(z-\xi)=\rho_1(z)\le r_1(z),$$
showing that, $z_0\in B_{\overline q}[z,r_1(z)].$

Also
$$
q(z_0-z)=\mu(z_0-z)=\mu(P(\xi-z))\le\mu(\xi-z)=\rho_2(z)\le r_2(z),$$
that is, $z_0\in B_{q}[z,r_2(z)].$

Hence,
$$
z_0\in \bigcap_{z\in Z}\left(B_{q}[z,r_2(z)]\cap B_{\bar q}[z,r_1(z)]\right),$$
in contradiction  to \eqref{eq.asym-ety}. By Proposition \ref{p1.ell-infty-inj}.5, this shows that the asymmetric normed space  $Z$ is not injective.
\hfill $\Box$\end{pf}

\end{document}